\documentclass{matharticle}
\usepackage{graphicx}
\usepackage{url}
\newtheorem{theorem}{Theorem}[section]
\newtheorem{lemma}[theorem]{Lemma}
\usepackage{hyperref}
\hypersetup{hidelinks}
\theoremstyle{definition}

\newtheorem{example}[theorem]{Example}

\newtheorem{corollary}[theorem]{Corollary}

\theoremstyle{remark}
\newtheorem{remark}[theorem]{Remark}

\numberwithin{equation}{section}

\begin{document}

\title[Counterexample to the strong union of balls conjecture]{A counterexample to the strong version of the union of uniform
closed balls conjecture in dimensions three and higher}

\author{Chadi Nour}

\address{Department of Computer Science and Mathematics, Lebanese American University, Byblos Campus, P.O. Box 36, Byblos, Lebanon}
\email{cnour@lau.edu.lb}



\begin{abstract}
Fifteen years after its formulation and following its recent proof in the
plane, we show that the strong version of the union of uniform closed balls
conjecture fails in every dimension $n\geq 3$. For every $n\geq 3$, $r>0$, and $\varepsilon>0$, we construct a compact set
$S\subset\mathbb R^n$ satisfying the interior $r$-sphere condition and
containing a point that belongs to no closed ball of radius
$(1/2+\varepsilon)r$ contained in $S$. Consequently, the known universal radius $r/2$ is optimal in every dimension
$n\geq 3$, even for compact sets, whereas the optimal radius in the plane is
$r/\sqrt{3}$. The construction uses simplices with coincident circumcenter and incenter and
yields an exact formula for the largest radius of a contained ball covering
the center; the regular simplex recovers the original example motivating
the conjecture.
\end{abstract}
\maketitle

\section{Introduction and main result}

Let $S\subset\mathbb R^n$ be a nonempty closed set and let $r>0$.
We recall that $S$ is said to satisfy the \emph{interior $r$-sphere
condition} if every point of its boundary lies on the boundary of a
closed ball of radius $r$ contained in $S$. The interior sphere
condition is a classical geometric hypothesis in elliptic boundary-point
estimates \cite{Hopf1952}, the regularity theory of attainable sets and
minimum-time functions \cite{CF2006}, and perimeter estimates for
reachable sets \cite{CC2006}.

One can easily see that if $S$ is the union of closed balls of radius
$r$, that is,
\[
S=S_r:=\bigcup\bigl\{\overline B(c;r):
                     \overline B(c;r)\subset S\bigr\},
\]
where $\overline B(c;r)$ denotes the closed ball centered at $c$ with
radius $r$, then $S$ satisfies the interior $r$-sphere condition.
The converse is not necessarily true. Indeed, in their study of
proximal smoothness and the exterior sphere condition, Nour, Stern,
and Takche gave a planar set satisfying the interior unit-sphere
condition whose origin belongs to no closed unit ball contained in
the set \cite[Example~4.1]{NST2009JCA}; see also
\cite[Example~6]{NST2009CC}. In 2009, they formulated the conjecture that the interior $r$-sphere condition nevertheless implies
$S=S_{r'}$ for some $r'>0$
\cite[Conjecture~4.4]{NST2009JCA}, also stated in
\cite[Conjecture~1]{NST2009CC}. This is now called the
\emph{weak version} of the conjecture. They also proved it when $S$
is wedged, or epi-Lipschitz, and has compact boundary
\cite[Corollary~4.2]{NST2009JCA}.

In 2011, Nour, Stern, and Takche proved the weak conjecture without
these additional assumptions \cite{NST2011JCA}. More precisely,
they showed that every nonempty closed set satisfying the interior
$r$-sphere condition satisfies
\begin{equation}\label{eq:weak-history}
S=S_{r/2}.
\end{equation}
Their proof uses nonsmooth and proximal analysis and also leads to
a variable-radius extension. For background on nonsmooth and variational analysis, we refer
to the monographs \cite{CLSW1998,Mordukhovich2006I,Mordukhovich2006II,RW1998,Thibault2023}.

In the same year, they generalized the original planar example to
every dimension $n\geq2$ \cite[Example~2]{NST2011DCDS}.
Their construction uses a regular simplex and gives a set satisfying
the interior $r$-sphere condition for which the origin belongs to no
contained closed ball of radius greater than
\begin{equation}\label{eq:target}
\rho_n(r)=\frac{nr}{2\sqrt{n^2-1}}.
\end{equation}
This example motivated the \emph{strong version} of the conjecture,
formulated in \cite[Conjecture~3]{NST2011DCDS} and also stated in
\cite[Conjecture~2.5]{NST2011JCA}: every nonempty closed set satisfying
the interior $r$-sphere condition should satisfy
$S=S_{\rho_n(r)}$. Since $\rho_n(r)\downarrow r/2$ as $n\to\infty$, the regular-simplex
examples already show that $r/2$ is the largest radius that can work
uniformly over all dimensions. They do not, however, determine
whether the radius can be increased when the dimension is fixed.

More recently, Nour and Takche proved the strong conjecture in the
plane \cite[Theorem~1.4]{NTplane}. They showed that every nonempty
closed set $S\subset\mathbb R^2$ satisfying the interior $r$-sphere
condition is the union of contained closed balls of radius
$r/\sqrt3$. The original planar example shows that this radius is
optimal.

The purpose of the present paper is to show that the situation
changes in all dimensions three and higher. We disprove the strong
conjecture and prove that the radius $r/2$ in
\eqref{eq:weak-history} is optimal in each fixed dimension $n\geq3$,
even when $S$ is required to be compact. To state the result, define
\[
R_{\max}(S,x)
=
\sup\bigl\{a>0:
  \text{there exists }c\in\mathbb R^n
  \text{ with }x\in\overline B(c;a)\subset S\bigr\}.
\]

\begin{theorem}\label{thm:main}
Let $n\geq 3$, $r>0$, and $\varepsilon>0$. There exists a nonempty compact set $S\subset\mathbb{R}^n$
containing the origin and satisfying the interior
$r$-sphere condition such that
\[
R_{\max}(S,0)<\left(\frac12+\varepsilon\right)r.
\]
In particular, the strong conjecture is false in every dimension
$n\geq3$, and the universal radius $r/2$ cannot be increased,
even for compact sets.
\end{theorem}

The construction extends the geometry of
\cite[Example~2]{NST2011DCDS} to simplices whose circumcenter and
incenter coincide, without requiring regularity. For a simplex of
unit circumradius and inradius $h$, we construct a compact set
$S_\Delta$ satisfying the interior unit-sphere condition and prove
the exact identity
\[
R_{\max}(S_\Delta,0)=\frac{1}{2\sqrt{1-h^2}}.
\]
The regular simplex has $h=1/n$ and recovers the original example.
In dimension three, nonregular tetrahedra of this type have
$h<1/3$ and yield counterexamples consisting of one small ball
and four unit balls. The optimal covering radii of these sets tend to $1/2$ from above
as the tetrahedra degenerate. To obtain compact counterexamples
in every dimension $n\geq3$, we construct an explicit cyclic
deformation of the regular simplex that preserves the coincidence
of the circumcenter and incenter while allowing the inradius to
tend to zero. Applying the simplex construction to this family
proves Theorem~\ref{thm:main}. Such a degeneration is impossible
in the plane. Indeed, a triangle whose circumcenter and incenter coincide
must be equilateral. Noncompact closed counterexamples in higher
dimensions can also be obtained by taking Cartesian products of
suitable three-dimensional examples with Euclidean spaces; see
Remark~\ref{products}.

Note that we normalize $r=1$ in the proofs. The general case follows by
dilation, which transforms the interior unit-sphere condition
into the interior $r$-sphere condition and multiplies every
covering radius by $r$.

The paper is organized as follows. Section~\ref{simplex} presents
the simplex construction, verifies the interior sphere condition,
and determines the exact optimal covering radius. It also includes
an explicit three-dimensional counterexample illustrating the construction.
Section~\ref{higher} constructs degenerating nonregular simplices
in every dimension $n\geq3$ and proves Theorem~\ref{thm:main}.

\section{The simplex construction and its optimal covering radius}
\label{simplex}

Throughout, $\langle\cdot,\cdot\rangle$ and $\|\cdot\|$ denote the
Euclidean inner product and norm, respectively. We write
$\mathbb{S}^{n-1}$ for the unit sphere in $\mathbb{R}^n$ and
$\operatorname{co} A$ for the convex hull of a set $A$.
We write $B(x;a)$ for the open ball centered at $x$ with radius $a$,
and $\mathbb{S}(x;a)$ for the sphere centered at $x$ with radius $a$.

Let $n\geq2$, and consider a nondegenerate simplex
\[
\Delta=\operatorname{co}\{u_0,\ldots,u_n\}\subset\mathbb R^n,
\]
where $u_0,\ldots,u_n\in\mathbb S^{n-1}$ and
$0\in\operatorname{int}\Delta$. Assume that all facets of $\Delta$
lie at the same distance $h>0$ from the origin. Thus the
circumcenter and incenter of $\Delta$ coincide at the origin,
and its circumradius and inradius are $1$ and $h$, respectively. For $j\in\{0,\ldots,n\}$, let $F_j$ be the facet opposite $u_j$,
and let $n_j$ be its outward unit normal. Then
\begin{equation}\label{eq:facets}
\langle u_i,n_j\rangle=h \quad(i\neq j),
\;\;\;
\langle u_j,n_j\rangle<h,
\end{equation}
and
\[
\Delta=\bigcap_{j=0}^n
\{x\in\mathbb R^n:\langle x,n_j\rangle\leq h\}.
\]
In particular, $\overline B(0;h)\subset\Delta$. 

Recall that a simplex is called \emph{regular} if all its edges
have the same length.

\begin{lemma}\label{lem:euler}
One has $h\leq1/n$, with equality if and only if $\Delta$ is regular.
Moreover, for every $z\in\mathbb R^n$,
\begin{equation}\label{eq:support}
\max_{0\leq i\leq n}\langle z,u_i\rangle\geq h\|z\|.
\end{equation}
\end{lemma}

\begin{proof}
Since $0\in\operatorname{int}\Delta$, there exist
$\lambda_0,\ldots,\lambda_n>0$ such that
\[
\sum_{i=0}^n\lambda_i=1
\;\;\hbox{and}\;\;
\sum_{i=0}^n\lambda_i u_i=0.
\]
Taking the inner product with $n_j$ and using
\eqref{eq:facets}, we obtain
\[
0=\lambda_j\langle u_j,n_j\rangle+(1-\lambda_j)h.
\]
Because $u_j$ and $n_j$ are unit vectors,
\[
-1\leq\langle u_j,n_j\rangle
=-h\frac{1-\lambda_j}{\lambda_j},
\]
and hence
\[
\lambda_j\geq\frac{h}{1+h}.
\]
Summing over $j$ gives
\[
1\geq\frac{(n+1)h}{1+h},
\]
which is equivalent to $h\leq1/n$. If $h=1/n$, all the preceding inequalities are equalities.
Consequently,
\[
\lambda_j=\frac1{n+1}
\;\;\text{and}\;\;
\langle u_j,n_j\rangle=-1,
\]
so $n_j=-u_j$. Equation \eqref{eq:facets} then gives
$\langle u_i,u_j\rangle=-1/n$ for $i\neq j$, and $\Delta$ is
regular. Conversely, a regular simplex of unit circumradius has
$n_j=-u_j$ and $h=1/n$. Finally, if $z\neq0$, then
$hz/\|z\|\in\overline B(0;h)\subset\Delta$. Since a linear
functional attains its maximum over a simplex at a vertex,
\[
h\|z\|
=\left\langle z,\frac{hz}{\|z\|}\right\rangle
\leq\max_{0\leq i\leq n}\langle z,u_i\rangle.
\]
The case $z=0$ is immediate.
\end{proof}

Define
\begin{equation}\label{eq:construction}
\rho_0:=\frac{h}{\sqrt{1-h^2}},
\;\;\;
p_i:=\rho_0u_i,
\;\;\;
c_j:=\frac{n_j}{\sqrt{1-h^2}}
     =\frac{\rho_0}{h}n_j,
\end{equation}
and set
\[
S_\Delta
:=
\overline B(0;\rho_0)
\cup
\bigcup_{j=0}^n\overline B(c_j;1).
\]
These quantities are well defined because
$0<h\leq1/n\leq1/2$.

\begin{lemma}\label{lem:contacts}
For every $j$,
\[
\|c_j\|^2=1+\rho_0^2.
\]
For $i\neq j$,
\[
\langle c_j,p_i\rangle=\rho_0^2,
\;\;\;
\|p_i-c_j\|=1,
\;\;\;
\langle c_j-p_i,p_i\rangle=0.
\]
Furthermore, $\|p_i-c_i\|>1$ for every $i$.
\end{lemma}

\begin{proof}
Since $\|n_j\|=1$,
\[
\|c_j\|^2=\frac1{1-h^2}=1+\rho_0^2.
\]
For $i\neq j$, equation \eqref{eq:facets} gives
\[
\langle c_j,p_i\rangle
=\frac{\rho_0^2}{h}\langle n_j,u_i\rangle
=\rho_0^2.
\]
Therefore
\[
\|p_i-c_j\|^2
=\rho_0^2+(1+\rho_0^2)-2\rho_0^2
=1,
\]
and
\[
\langle c_j-p_i,p_i\rangle
=\rho_0^2-\|p_i\|^2=0.
\]
For every $i$, the strict inequality in
\eqref{eq:facets} yields
\[
\|p_i-c_i\|^2
=
1+2\rho_0^2
-\frac{2\rho_0^2}{h}\langle n_i,u_i\rangle
>1.
\]
\end{proof}

\begin{theorem}\label{thm:exact}
The set $S_\Delta$ is compact and satisfies the interior
unit-sphere condition. Moreover,
\begin{equation}\label{eq:exact}
R_{\max}(S_\Delta,0)
=
\frac12\sqrt{1+\rho_0^2}
=
\frac1{2\sqrt{1-h^2}}.
\end{equation}
The supremum is attained by each of the balls
\[
\overline B\left(\frac{c_j}{2};\frac{\|c_j\|}{2}\right),
\;\;\; j=0,\ldots,n.
\]
\end{theorem}

\begin{proof}
Compactness follows because $S_\Delta$ is a finite union of
closed balls. We first verify the interior unit-sphere condition.
Let $w\in\mathbb S^{n-1}$. Since $0\in\operatorname{int}\Delta$,
the ray $\{tw:t\geq0\}$ meets the boundary of $\Delta$ at a point
$y$ belonging to some facet $F_j$. As
$\Delta\subset\overline B(0;1)$, we have $\|y\|\leq1$, and hence
\[
\langle w,n_j\rangle
=\frac{\langle y,n_j\rangle}{\|y\|}
=\frac{h}{\|y\|}
\geq h.
\]
Consequently, for $q=\rho_0w$,
\[
\|q-c_j\|^2
=
\rho_0^2+\frac{\rho_0^2}{h^2}
-\frac{2\rho_0^2}{h}\langle w,n_j\rangle\\
\leq
\frac{\rho_0^2}{h^2}-\rho_0^2
=1.
\]
Thus every point of $\mathbb S(0;\rho_0)$ belongs to at least
one of the displayed closed unit balls. Every boundary point of $S_\Delta$ lies on the boundary of
at least one of its generating balls. If it lies on a unit
sphere, the required condition follows from the corresponding unit ball. Otherwise, it lies on $\mathbb S(0;\rho_0)$ and
therefore belongs to a contained closed unit ball by the
preceding argument. Since it is a boundary point of $S_\Delta$,
it cannot be interior to that unit ball and thus lies on its
boundary. This proves the sphere condition at every boundary
point. We next show that each $p_i$ is a boundary point of $S_\Delta$.
For $\epsilon>0$, the point $(1+\epsilon)p_i$ lies outside
$\overline B(0;\rho_0)$. For $j\neq i$, Lemma~\ref{lem:contacts}
gives
\[
\|(1+\epsilon)p_i-c_j\|^2
=
\|p_i-c_j\|^2
+2\epsilon\langle p_i-c_j,p_i\rangle
+\epsilon^2\rho_0^2\\
=
1+\epsilon^2\rho_0^2
>1.
\]
Since $\|p_i-c_i\|>1$, the point $(1+\epsilon)p_i$ also lies
outside $\overline B(c_i;1)$ for all sufficiently small
$\epsilon>0$. Hence points outside $S_\Delta$ approach $p_i$,
and $p_i\in\operatorname{bdry}S_\Delta$. We now prove the upper bound in \eqref{eq:exact}. Put
\[
M:=\frac1{2\sqrt{1-h^2}}=\frac{\rho_0}{2h},
\]
and suppose
\[
0\in\overline B(z;T)\subset S_\Delta.
\]
Writing $t=\|z\|$, we have $t\leq T$. Since each $p_i$ is a
boundary point of $S_\Delta$, it cannot belong to
$B(z;T)\subset\operatorname{int}S_\Delta$. Therefore
\[
T^2\leq\|z-p_i\|^2
\;\;\text{for every }i.
\]
Choosing an index that maximizes $\langle z,u_i\rangle$ and
using \eqref{eq:support}, we obtain
\begin{equation}\label{eq:radius-upper}
T^2
\leq
\rho_0^2+t^2
-2\rho_0\max_i\langle z,u_i\rangle
\leq
\rho_0^2+t^2-2\rho_0ht.
\end{equation}
Together with $t^2\leq T^2$, this implies
\[
t\leq\frac{\rho_0}{2h}=M.
\]
The function $t\mapsto
f(t)=\rho_0^2+t^2-2\rho_0ht$ is convex on $[0,M]$, with
\[
f(0)=\rho_0^2
\;\;\text{and}\;\;
f(M)=\rho_0^2+M^2-2\rho_0hM=M^2.
\]
Because $h\leq1/2$, we have $\rho_0\leq M$. Thus both endpoint
values are at most $M^2$, and \eqref{eq:radius-upper} yields
$T\leq M$. Finally, fix $j$ and consider the ball
\[
\beta_j:=
\overline B\left(\frac{c_j}{2};\frac{\|c_j\|}{2}\right).
\]
For every $x\in\mathbb R^n$, the identity
$\|c_j\|^2=1+\rho_0^2$ gives
\begin{align*}
\frac12\bigl(\|x\|^2-\rho_0^2\bigr)
+\frac12\bigl(\|x-c_j\|^2-1\bigr)
&=
\left\|x-\frac{c_j}{2}\right\|^2
+\frac{\|c_j\|^2}{4}
-\frac{1+\rho_0^2}{2}\\
&=
\left\|x-\frac{c_j}{2}\right\|^2
-\frac{\|c_j\|^2}{4}.
\end{align*}
For $x\in\beta_j$, the right-hand side is nonpositive.
Hence at least one of the two terms averaged on the left
is nonpositive, proving
\[
\beta_j
\subset
\overline B(0;\rho_0)\cup\overline B(c_j;1)
\subset S_\Delta.
\]
The ball $\beta_j$ contains the origin on its boundary and
has radius $\|c_j\|/2=M$. This proves attainment and
\eqref{eq:exact}.
\end{proof}

\begin{corollary}\label{cor:regular}
One has
\[
R_{\max}(S_\Delta,0)\leq
\rho_n(1)=\frac{n}{2\sqrt{n^2-1}},
\]
with equality if and only if $\Delta$ is regular.
In dimension two, $\Delta$ is necessarily regular.
Consequently, if $n\geq3$ and $\Delta$ is nonregular, then
$S_\Delta$ is a counterexample to the strong conjecture.
\end{corollary}

\begin{proof}
The function $h\mapsto1/(2\sqrt{1-h^2})$ is strictly increasing
on $(0,1)$. The inequality and its equality case therefore follow
from Lemma~\ref{lem:euler} and Theorem~\ref{thm:exact}. When $n=2$, the three sides of $\Delta$ are chords of the unit
circle whose supporting lines lie at the same distance $h$
from the origin. Each side has length $2\sqrt{1-h^2}$,
so $\Delta$ is equilateral.
\end{proof}

\begin{remark}
For a regular simplex, $h=1/n$ and $n_j=-u_j$. Thus
\[
\rho_0=\frac1{\sqrt{n^2-1}}
\;\;\text{and}\;\;
c_j=-\frac{n}{\sqrt{n^2-1}}u_j=-np_j.
\]
In this case, $S_\Delta$ is precisely the original example
of Nour, Stern, and Takche \cite[Example~2]{NST2011DCDS}.
\end{remark}

We next exhibit a nonregular tetrahedron whose circumcenter and
incenter coincide. The associated set provides an explicit compact
counterexample to the strong conjecture in $\mathbb R^3$.

\begin{example}\label{ex:exact}
Let
\[
\mathcal E
=\{(1,1,1),(1,-1,-1),(-1,1,-1),(-1,-1,1)\},
\]
and, for $s=(s_1,s_2,s_3)\in\mathcal E$, define
\[
u_s=\frac13(2s_1,2s_2,s_3),
\;\;\;
p_s=\frac{\sqrt2}{15}(2s_1,2s_2,s_3),\;\;\text{and}\;\;
c_s=-\frac{3\sqrt2}{10}(s_1,s_2,2s_3).
\]
Then
\[
S_\Delta
=
\overline B(0;\sqrt2/5)
\cup
\bigcup_{s\in\mathcal E}\overline B(c_s;1)
\]
is a compact set satisfying the interior unit-sphere condition, and
\[
R_{\max}(S_\Delta,0)
=\frac{3\sqrt3}{10}
<\frac3{4\sqrt2}
=\rho_3(1).
\]
Thus $S_\Delta$ is a counterexample to the strong conjecture. The construction is illustrated in panels~(a) and~(b) of
Figure~\ref{fig:counterexample}.
\begin{figure}[tb]
\centering
\includegraphics[ width=\textwidth,trim=50 6 50 6,clip]{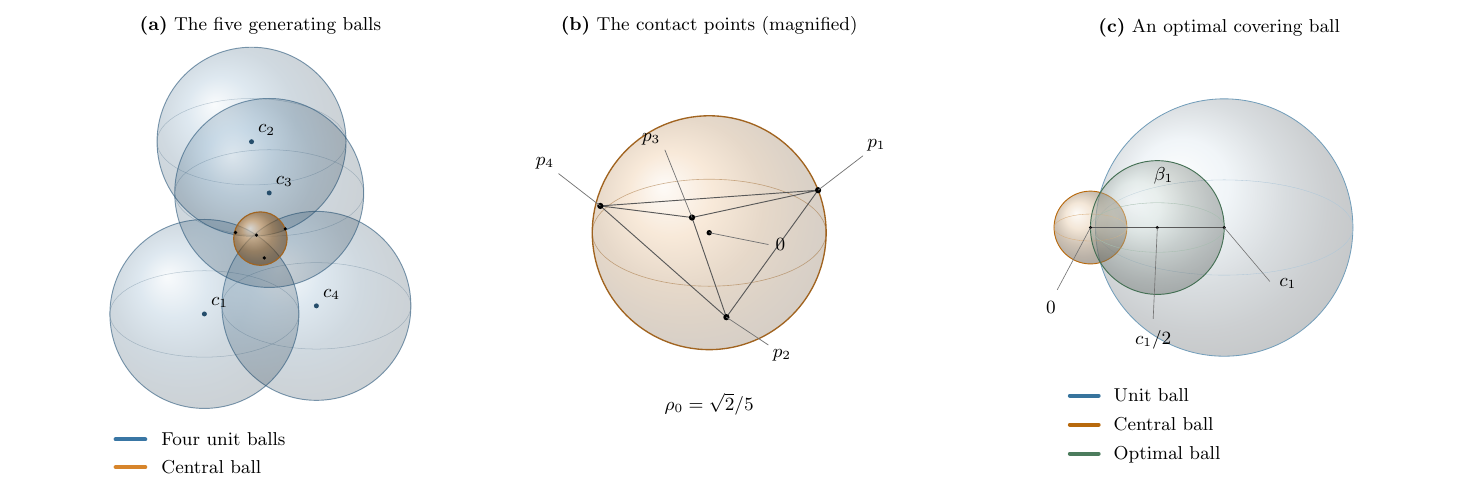}
\caption{The set $S_\Delta$ of Example~\ref{ex:exact}. Panel~(a) shows the four unit balls and the central ball as transparent
surfaces. Panel~(b) shows the boundary contact points and their
nonregular tetrahedron at a larger scale. Panel~(c) shows an optimal
covering ball $\beta_1=\overline{B}(c_1/2;3\sqrt{3}/10)$, contained
in the union of the central ball and $\overline{B}(c_1;1)$, with
the origin on its boundary. The indices $1,\ldots,4$ correspond
to the displayed order of the sign vectors in $\mathcal E$.
Hidden geometry is displayed for clarity.}
\label{fig:counterexample}
\end{figure}
To verify these assertions, first observe that the vectors $u_s$
are unit vectors and sum to zero. Enumerating them in the displayed
order, we have
\[
\det(u_2-u_1,u_3-u_1,u_4-u_1)
=-\frac{64}{27}\neq0.
\] They therefore form a nondegenerate tetrahedron $\Delta$ with
$0\in\operatorname{int}\Delta$. For each $s\in\mathcal E$, put
\[
n_s=-\frac1{\sqrt6}(s_1,s_2,2s_3).
\]
These are unit vectors and satisfy
\[
\langle u_t,n_s\rangle
=
\begin{cases}
  2/(3\sqrt6),&t\neq s,\\
  -2/\sqrt6,&t=s.
\end{cases}
\]
Hence $n_s$ is the outward unit normal to the facet opposite $u_s$,
and all four facets lie at distance $h=2/(3\sqrt6)$ from the origin.
The circumcenter and incenter therefore coincide at zero.
In particular,
\[
h^2=\frac2{27},
\;\;\;
\rho_0^2=\frac{h^2}{1-h^2}=\frac2{25},\;\;\text{and}\;\;
\frac{n_s}{\sqrt{1-h^2}}=c_s.
\]
Thus this is precisely the construction of this section. Theorem~\ref{thm:exact} proves the interior unit-sphere condition and gives
\[
R_{\max}(S_\Delta,0)^2
=\frac1{4(1-h^2)}
=\frac{27}{100}.
\]
An optimal covering ball is shown in panel~(c) of Figure~\ref{fig:counterexample}. Finally,
\[
\rho_3(1)^2-R_{\max}(S_\Delta,0)^2
=\frac9{32}-\frac{27}{100}
=\frac9{800}>0,
\]
which proves the strict inequality. \end{example}

\section{Compact counterexamples in higher dimensions}\label{higher}

We now construct simplices with coincident circumcenter and incenter
whose inradii tend to zero in every fixed dimension $n\geq3$.
Applying the construction of Section~\ref{simplex} to these simplices
will prove Theorem~\ref{thm:main}.

Fix $n\geq3$ and put $m=n+1$. Consider the $n$-dimensional Euclidean
space
\[
H=\left\{x\in\mathbb R^m:\sum_{i=0}^{m-1}x_i=0\right\},
\]
which we identify isometrically with $\mathbb R^n$.
Let $e_0,\ldots,e_{m-1}$ be the standard basis of $\mathbb R^m$,
write $\mathbf1=(1,\ldots,1)$, and set
\[
v_i=e_i-\frac{\mathbf1}{m},
\;\;\; i=0,\ldots,m-1.
\]
The vectors $v_i$ are the vertices of a regular simplex in $H$,
centered at zero, and satisfy
\begin{equation}\label{eq:standard-simplex}
\sum_{i=0}^{m-1}v_i=0,
\;\;\;
\|v_i\|^2=\frac{n}{m},
\;\;\text{and}\;\;
\langle v_i,v_j\rangle=-\frac1m
\quad(i\neq j).
\end{equation}
Define $q,s\in\mathbb R^m$ by
\[
q_i=\cos\left(\frac{2\pi i}{m}\right)
\;\;\text{and}\;\;
s_i=\sin\left(\frac{2\pi i}{m}\right),
\;\;\; i=0,\ldots,m-1.
\]
The elementary trigonometric sums give
\[
\sum_{i=0}^{m-1}q_i=\sum_{i=0}^{m-1}s_i=0,
\;\;\;
\langle q,s\rangle=0,
\;\;\text{and}\;\;
\|q\|^2=\|s\|^2=\frac m2.
\]
Indeed, these identities follow by taking real and imaginary
parts of the geometric sums
\[
\sum_{j=0}^{m-1}
\exp\left(\frac{2\pi kj}{m}\mathrm{i}\right)=0,
\;\; k=1,2,
\]
which vanish because $\exp(2\pi k\mathrm{i}/m)\ne1$ for
$k=1,2$, as $m\geq4$, and using the double-angle formulas.
Since $q,s\in H$ are nonzero and orthogonal,
$W:=\operatorname{span}\{q,s\}$ is a two-dimensional subspace of $H$. Let $P\colon H\to W$ be the orthogonal projection.
Then
\[
Px=\frac2m\bigl(\langle x,q\rangle q
                 +\langle x,s\rangle s\bigr),
\;\;\; x\in H.
\]
For $0<t\leq1$, define
\[
T_t=P+t(I_H-P),
\]
where $I_H$ is the identity on $H$. The map $T_t$ is self-adjoint
and invertible, with
\[
T_t^{-1}=P+t^{-1}(I_H-P).
\]
It fixes $W$ and contracts its orthogonal complement in $H$ by
the factor $t$. The assumption $n\geq3$ ensures that the orthogonal complement
of $W$ in $H$ has dimension $n-2\geq1$, so this deformation
is nontrivial for $0<t<1$. In contrast, when $n=2$, we have
$W=H$ and $T_t=I_H$ for every $t>0$.

\begin{lemma}\label{lem:cyclic}
For $0<t\leq1$, put
\[
N_t=\sqrt{\frac{2+(n-2)t^2}{m}}
\;\;\;\text{and}\;\;\;
u_i(t)=\frac{T_tv_i}{N_t},
\;\;\; i=0,\ldots,m-1.
\]
Then
\[
\Delta_t=\operatorname{co}\{u_0(t),\ldots,u_{m-1}(t)\}
\]
is a nondegenerate simplex of unit circumradius whose circumcenter
and incenter coincide at zero. Its inradius is
\begin{equation}\label{eq:cyclic-inradius}
h_t=
\frac1{\sqrt{(2+(n-2)t^2)(2+(n-2)t^{-2})}}.
\end{equation}
Moreover,
\[
h_1=\frac1n,
\;\;\;
h_t<\frac1n \quad(0<t<1),
\;\;\text{and}\;\;
\lim_{t\downarrow0}h_t=0.
\]
In particular, $\Delta_t$ is nonregular whenever $0<t<1$.
\end{lemma}

\begin{proof}
Since $q,s\in H$, we have
\[
\langle v_i,q\rangle=q_i
\;\;\text{and}\;\;
\langle v_i,s\rangle=s_i.
\]
Thus
\[
Pv_i=\frac2m(q_iq+s_is),
\]
and the orthogonality and norm identities for $q$ and $s$ yield
\[
\|Pv_i\|^2
=\frac2m(q_i^2+s_i^2)
=\frac2m.
\]
Using \eqref{eq:standard-simplex}, we obtain
\[
\|(I_H-P)v_i\|^2
=\|v_i\|^2-\|Pv_i\|^2
=\frac{n-2}{m}.
\]
The two components are orthogonal, so
\[
\|T_tv_i\|^2
=\|Pv_i\|^2+t^2\|(I_H-P)v_i\|^2
=\frac{2+(n-2)t^2}{m}
=N_t^2.
\]
Therefore every $u_i(t)$ is a unit vector. Since $T_t$ is invertible, it preserves affine independence.
The vectors $u_i(t)$ consequently form a nondegenerate simplex.
Moreover,
\[
\sum_{i=0}^{m-1}u_i(t)=0,
\]
so zero has equal positive barycentric coordinates and belongs
to $\operatorname{int}\Delta_t$. Since all vertices are unit
vectors, zero is its circumcenter. To determine the facet normals, put
\[
L_t=\sqrt{\frac{2+(n-2)t^{-2}}{m}}
\;\;\text{and}\;\;
n_j(t)=-\frac{T_t^{-1}v_j}{L_t}.
\]
The same orthogonal decomposition gives
\[
\|T_t^{-1}v_j\|^2
=\frac{2+(n-2)t^{-2}}{m}
=L_t^2,
\]
so $n_j(t)$ is a unit vector. By self-adjointness,
\[
\langle u_i(t),n_j(t)\rangle
=
-\frac{\langle T_tv_i,T_t^{-1}v_j\rangle}{N_tL_t}\\
=
-\frac{\langle v_i,v_j\rangle}{N_tL_t}.
\]
Hence
\[
\langle u_i(t),n_j(t)\rangle
=
\begin{cases}
  1/(mN_tL_t),&i\neq j,\\
  -n/(mN_tL_t),&i=j.
\end{cases}
\]
Since
\[
\langle u_j(t),n_j(t)\rangle
=-\frac{n}{mN_tL_t}<\frac{1}{mN_tL_t},
\]
the vertex $u_j(t)$ lies strictly on the inner side of the
hyperplane containing the opposite facet. Thus $n_j(t)$ is
the outward unit normal to that facet, and every facet lies
at the same positive distance
\[
h_t=\frac{1}{mN_tL_t}
\]
from zero. Since zero is interior to the simplex, it is also
the incenter. Substituting the expressions for $N_t$ and $L_t$
proves \eqref{eq:cyclic-inradius}. At $t=1$, equation \eqref{eq:cyclic-inradius} gives $h_1=1/n$.
For $0<t<1$, we have $t^2+t^{-2}>2$ and $n-2>0$, whence
\begin{eqnarray*}
(2+(n-2)t^2)(2+(n-2)t^{-2})&=&4+(n-2)^2+2(n-2)(t^2+t^{-2})\\&>&4+(n-2)^2+4(n-2)=n^2.
\end{eqnarray*}
Thus $h_t<1/n$, and Lemma~\ref{lem:euler} shows that
$\Delta_t$ is nonregular. Finally, the denominator in
\eqref{eq:cyclic-inradius} tends to infinity as $t\downarrow0$,
so $h_t\to0$.
\end{proof}

\begin{remark}
For $n=3$ and $t=1/2$, the simplex $\Delta_t$ is congruent to the tetrahedron
of Example~\ref{ex:exact}. Indeed, for $i\ne j$, a direct computation gives
\[
\langle u_i(1/2),u_j(1/2)\rangle=
\begin{cases}
-\dfrac79, & \text{if } |i-j|=2,\\[7pt]
-\dfrac19, & \text{if } |i-j| \text{ is odd}.
\end{cases}
\]
These are exactly the inner products $\langle u_s,u_{s'}\rangle$
of Example~\ref{ex:exact}, with $\{u_0,u_2\}$ and $\{u_1,u_3\}$
corresponding to the two pairs of sign vectors with equal third
coordinate. Equal Gram matrices give congruent simplices.
\end{remark}

Now we prove Theorem~\ref{thm:main}.

\begin{proof}[Proof of Theorem~\ref{thm:main}]
First take $r=1$. For each $0<t<1$, apply the construction
of Section~\ref{simplex} to $\Delta_t$, and denote the resulting
set by $S_{\Delta_t}$. More explicitly, with
\[
\rho_0(t)=\frac{h_t}{\sqrt{1-h_t^2}},
\;\;\;
c_j(t)=\frac{n_j(t)}{\sqrt{1-h_t^2}},
\]
we set
\[
S_{\Delta_t}=
\overline B_H(0;\rho_0(t))
\cup
\bigcup_{j=0}^{m-1}\overline B_H(c_j(t);1).
\]
Here the balls are taken in the Euclidean space $H$, identified
isometrically with $\mathbb R^n$. By Theorem~\ref{thm:exact}, $S_{\Delta_t}$ is compact, satisfies the interior
unit-sphere condition, and
\[
R_{\max}(S_{\Delta_t},0)
=\frac1{2\sqrt{1-h_t^2}}.
\]
Lemma~\ref{lem:cyclic} therefore gives
\[
\lim_{t\downarrow0}R_{\max}(S_{\Delta_t},0)=\frac12.
\]
Given $\varepsilon>0$, choose $t\in(0,1)$ sufficiently small that
\[
R_{\max}(S_{\Delta_t},0)<\frac12+\varepsilon.
\]
For arbitrary $r>0$, let $S=rS_{\Delta_t}$. Dilation preserves compactness,
transforms the interior unit-sphere condition into the interior
$r$-sphere condition, and gives
\[
R_{\max}(S,0)
=rR_{\max}(S_{\Delta_t},0)
<\left(\frac12+\varepsilon\right)r.
\]
This proves the asserted construction. Since $\rho_n(r)>r/2$, choosing
\[
0<\varepsilon<\rho_n(1)-\frac12
\]
gives a compact set satisfying the interior $r$-sphere condition
whose origin belongs to no contained closed ball of radius
$\rho_n(r)$. Thus the strong conjecture fails. Finally, for any proposed universal constant $c>1/2$, choose $0<\varepsilon<c-1/2$. The same construction gives a compact set
whose origin belongs to no contained closed ball of radius $cr$.
Consequently, the universal radius $r/2$ cannot be increased.
\end{proof}

\begin{remark}\label{products}
Cartesian products provide another extension of the three-dimensional
construction. If $K\subset\mathbb R^3$ satisfies the interior
unit-sphere condition, then $K\times\mathbb R^{n-3}$ satisfies
the same condition for $n>3$. Indeed,
\[
\operatorname{bdry}(K\times\mathbb R^{n-3})
=
\operatorname{bdry}K\times\mathbb R^{n-3}.
\]
A contained unit ball tangent to $K$ at $x$ lifts to a contained
$n$-dimensional unit ball tangent to the product at $(x,y)$,
by choosing its center with second coordinate $y$. Moreover, orthogonal projection of any contained ball covering
$(0,0)$ gives a ball of the same radius contained in $K$ and
covering zero. Conversely, every such ball in $K$ lifts to a
ball of the same radius in the product. Therefore
\[
R_{\max}(K\times\mathbb R^{n-3},(0,0))
=
R_{\max}(K,0).
\]
Choosing $K=S_{\Delta_t}$ from the construction above in dimension
three, with $0<t<1$ sufficiently small that
$R_{\max}(K,0)<\rho_n(1)$, yields a counterexample in dimension $n$.
These product sets are noncompact, whereas the construction
above gives compact counterexamples directly in every dimension.
\end{remark}

\begin{remark}\label{rem:optimal-radius}
Combining Theorem~\ref{thm:main} with the universal covering
theorem of Nour, Stern, and Takche \cite{NST2011JCA} and the
planar theorem of Nour and Takche \cite{NTplane}, we conclude that
the largest universal common radius is
\[
\begin{cases}
r/\sqrt3,&n=2,\\
r/2,&n\geq3.
\end{cases}
\]
The planar sharpness follows from the original example
\cite[Example~4.1]{NST2009JCA}, while the higher-dimensional
sharpness holds even for compact sets. For $n\geq3$, the sharp radius $r/2$ is never attained as a pointwise
optimal covering radius. Indeed, let $S\subset\mathbb R^n$ be a
nonempty proper closed set satisfying the interior $r$-sphere condition.
For $x\in\operatorname{int}S$, let $d>0$ be the distance from $x$ to
$\operatorname{bdry}S$. By \cite[Remark~3.1]{NT2012JOTA},
\[
R_{\max}(S,x)\geq\frac{\sqrt{r^2+d^2}}2>\frac r2.
\]
For $x\in\operatorname{bdry}S$, the interior $r$-sphere condition
gives $R_{\max}(S,x)\geq r$.
\end{remark}

\section*{Use of artificial intelligence tools}
This work was conceived and directed by the author. AI assistants
(Claude Opus 5.5, Anthropic; GPT-6.1 Sol, OpenAI) were used for numerical exploration,
checking computations, simplifying proofs, and assisting with the
preparation of the manuscript in \LaTeX. The author takes full responsibility for all results, proofs,
and other content of the paper.

\bibliographystyle{amsplain}

\begin{thebibliography}{10}

\bibitem{CC2006}
P. Cannarsa and P. Cardaliaguet,
\textit{Perimeter estimates for reachable sets of control systems},
J. Convex Anal. \textbf{13} (2006), no. 2, 253--267.

\bibitem{CF2006}
P. Cannarsa and H. Frankowska,
\textit{Interior sphere property of attainable sets and time optimal
control problems},
ESAIM Control Optim. Calc. Var. \textbf{12} (2006), no. 2, 350--370.

\bibitem{CLSW1998}
F. H. Clarke, Yu. S. Ledyaev, R. J. Stern, and P. R. Wolenski,
\textit{Nonsmooth analysis and control theory},
Graduate Texts in Mathematics, vol.~178,
Springer-Verlag, New York, 1998.

\bibitem{Hopf1952}
E. Hopf,
\textit{A remark on linear elliptic differential equations of second order},
Proc. Amer. Math. Soc. \textbf{3} (1952), no. 5, 791--793.

\bibitem{Mordukhovich2006I}
B. S. Mordukhovich,
\textit{Variational analysis and generalized differentiation.
I. Basic theory},
Grundlehren der mathematischen Wissenschaften, vol.~330,
Springer-Verlag, Berlin, 2006.

\bibitem{Mordukhovich2006II}
B. S. Mordukhovich,
\textit{Variational analysis and generalized differentiation.
II. Applications},
Grundlehren der mathematischen Wissenschaften, vol.~331,
Springer-Verlag, Berlin, 2006.

\bibitem{NST2009JCA}
C. Nour, R. J. Stern, and J. Takche,
\textit{Proximal smoothness and the exterior sphere condition},
J. Convex Anal. \textbf{16} (2009), no. 2, 501--514.

\bibitem{NST2009CC}
C. Nour, R. J. Stern, and J. Takche,
\textit{The union of uniform closed balls conjecture},
Control Cybernet. \textbf{38} (2009), no. 4B, 1525--1534.

\bibitem{NST2011DCDS}
C. Nour, R. J. Stern, and J. Takche,
\textit{Generalized exterior sphere conditions and
$\varphi$-convexity},
Discrete Contin. Dyn. Syst. \textbf{29} (2011), no. 2, 615--622.

\bibitem{NST2011JCA}
C. Nour, R. J. Stern, and J. Takche,
\textit{Validity of the union of uniform closed balls conjecture},
J. Convex Anal. \textbf{18} (2011), no. 2, 589--600.

\bibitem{NT2012JOTA}
C. Nour and J. Takche,
\textit{On the union of closed balls property},
J. Optim. Theory Appl. \textbf{155} (2012), 376--389.

\bibitem{NTplane}
C. Nour and J. Takche,
\textit{Validity of the strong version of the union of uniform
closed balls conjecture in the plane},
J. Convex Anal., to appear; \href{https://arxiv.org/abs/2603.07588}{arXiv:2603.07588}.

\bibitem{RW1998}
R. T. Rockafellar and R. J.-B. Wets,
\textit{Variational analysis},
Grundlehren der mathematischen Wissenschaften, vol.~317,
Springer-Verlag, Berlin, 1998.

\bibitem{Thibault2023}
L. Thibault,
\textit{Unilateral variational analysis in Banach spaces},
Parts I and II, World Scientific, Singapore, 2023.

\end{thebibliography}

\end{document}